\documentclass[12pt,thmsa]{article}%
\usepackage{amsmath}
\usepackage{amsfonts}
\usepackage{amssymb}
\usepackage{sw20bams}
\usepackage{graphicx}%
\providecommand{\U}[1]{\protect\rule{.1in}{.1in}}
\begin{document}

\author{Steven Finch}
\title{Appell $F_{1}$ and Conformal Mapping }
\date{August 5, 2014}
\maketitle

\begin{abstract}
This is the last of a trilogy of papers on triangle centers. \ A fairly
obscure \textquotedblleft conformal center of gravity\textquotedblright\ is
computed for the class of all isosceles triangles. \ This calculation appears
to be new. \ A byproduct is the logarithmic capacity or transfinite diameter
of such, yielding results consistent with Haegi (1951).

\end{abstract}

\footnotetext{Copyright \copyright \ 2014 by Steven R. Finch. All rights
reserved.}Before discussing triangles, let us give both a review of \cite{Fi1}
and a preview involving a simpler region in the plane. \ 

Let $\Omega=\left\{  x+iy\in\mathbb{C}:y>0\text{, }x^{2}+y^{2}<1\right\}  $,
the interior of the upper half-disk of unit radius. \ Let $\Delta$ denote the
(full) disk of unit radius and $\Sigma$ denote the infinite horizontal strip
of width $\pi$. \ Define a function $\ell:\Omega\rightarrow\Sigma$ by
\cite{Kv}%
\[
\ell(z)=\ln\left(  \frac{(1+z)^{2}}{(1-z)^{2}}\right)  .
\]
The conformal map $f_{w}:\Omega\rightarrow\Delta$ given by
\[
f_{w}(z)=\frac{\exp(\ell(z))-\exp(\ell(w))}{\exp(\ell(z))-\exp(\overline
{\ell(w)})}%
\]
satisfies $f_{w}(w)=0$; it is well-known that $\ln\left\vert f_{w}%
(z)\right\vert $ is Green's function for $\Omega$. \ We deduce that%
\[
h(w)=\lim_{z\rightarrow w}\left\vert \frac{f_{w}(z)}{z-w}\right\vert
=\left\vert \frac{\exp(\ell(w))\ell^{\prime}(w)}{\exp(\ell(w))-\exp
(\overline{\ell(w)})}\right\vert
\]
where $\ell^{\prime}$ denotes the derivative of $\ell$. \ Restricting
attention to the $y$-axis only, we have%
\[
h(iy)=\left\vert \frac{\frac{2\left(  1+iy\right)  }{\left(  1-iy\right)
^{2}}+\frac{2\left(  1+iy\right)  ^{2}}{\left(  1-iy\right)  ^{3}}}%
{\frac{\left(  1+iy\right)  ^{2}}{\left(  1-iy\right)  ^{2}}-\frac{\left(
i+y\right)  ^{2}}{\left(  -i+y\right)  ^{2}}}\right\vert =\left\vert
\frac{\left(  1+iy\right)  ^{2}}{2y\left(  1-y^{2}\right)  }\right\vert
=\frac{1+y^{2}}{2y\left(  1-y^{2}\right)  }%
\]
for $0<y<1$. \ Minimizing this expression, it follows that%
\[
iy_{0}=i\sqrt{-2+\sqrt{5}}=(0.4858682717566456781828638...)i
\]
is the \textbf{inner conformal center} (what was called the least capacity
point in \cite{Fi1}) for $\Omega$. Also,%
\[
\frac{1}{h(iy_{0})}=\sqrt{-22+10\sqrt{5}}=0.6005662120015552157733894...
\]
is the \textbf{maximum inner radius} of $\Omega$ \cite{PS0, Ha}. This
concludes our review.

Let $\Omega^{c}$ denote the complement of the closure of $\Omega$. \ Inverting
a function \cite{Pm1}%
\[
\frac{1}{z}=\frac{\left(  \frac{\sqrt{3}}{2}-\frac{i}{2}\right)  \left(
1+\frac{1}{w}\right)  ^{2/3}+\left(  \frac{\sqrt{3}}{2}+\frac{i}{2}\right)
\left(  1-\frac{1}{w}\right)  ^{2/3}}{\left(  1+\frac{1}{w}\right)
^{2/3}-\left(  1-\frac{1}{w}\right)  ^{2/3}}%
\]
in terms of $w$ yields%
\[
\frac{\left(  1+\frac{1}{w}\right)  ^{2/3}}{\left(  1-\frac{1}{w}\right)
^{2/3}}=\frac{\frac{\sqrt{3}}{2}+\frac{i}{2}+\frac{1}{z}}{-\frac{\sqrt{3}}%
{2}+\frac{i}{2}+\frac{1}{z}}=m(z).
\]
The conformal map $g:\Delta\rightarrow\Omega^{c}$
\[
g(z)=\frac{1+2m(z)^{3/2}+m(z)^{3}}{-1+m(z)^{3}}=\frac{4}{3\sqrt{3}}\frac{1}%
{z}+\frac{2i}{3\sqrt{3}}+O\left(  z\right)
\]
satisfies $g(0)=\infty$, has positive leading Laurent coefficient, and is
unique in this regard. \ The constant term of the series expansion%
\[
\frac{2i}{3\sqrt{3}}=(0.3849001794597505096727658...)i
\]
is the \textbf{outer conformal center}, which clearly lies in $\Omega$ but is
not the same as the inner conformal center. \ A motivating feature is%
\[
\frac{1}{2\pi}%
{\displaystyle\int\limits_{0}^{2\pi}}
g\left(  (1-\varepsilon)e^{it}\right)  dt=\frac{2i}{3\sqrt{3}}%
\]
for $\varepsilon>0$, but the literature is small \cite{PS1, PS2, Pm2, BWY,
Ct}. \ The leading coefficient (of $1/z$) is the \textbf{outer radius} of
$\Omega$ \cite{PS0, Lp}:%
\[
\frac{4}{3\sqrt{3}}=0.7698003589195010193455317...
\]
which is also known as the \textbf{logarithmic capacity }or
\textbf{transfinite diameter} of $\Omega$. \ An alternative definition is
\cite{Ra1}
\[
\lim_{n\rightarrow\infty}\,\max_{z_{1},z_{2},...,z_{n}\in\Omega}\,\left(
{\displaystyle\prod\limits_{j<k}}
|z_{j}-z_{k}|\right)  ^{\frac{2}{n(n-1)}},
\]
that is, the maximal geometric mean of pairwise distances for $n$ points in
$\Omega$, in the limit as $n\rightarrow\infty$. \ This constitutes a
fascinating collision of ideas from potential theory; replacing a geometric
mean by an arithmetic mean seems to be an open computational issue.

As a conclusion to our preview, the Appell $F_{1}$ function can be written as
a definite integral \cite{By, MSR}%
\[
F_{1}(a,b,b^{\prime},c;x,y)=\frac{\Gamma(c)}{\Gamma(a)\Gamma(c-a)}%
{\displaystyle\int\limits_{0}^{1}}
s^{a-1}(1-s)^{c-a-1}(1-sx)^{-b}(1-sy)^{-b^{\prime}}ds,
\]%
\[%
\begin{array}
[c]{lllll}%
\left\vert x\right\vert <1, &  & \left\vert y\right\vert <1, &  &
\operatorname{Re}(c)>\operatorname{Re}(a)>0
\end{array}
\]
as well as a double hypergeometric series
\[
F_{1}(a,b,b^{\prime},c;x,y)=%
{\displaystyle\sum\limits_{m=0}^{\infty}}
{\displaystyle\sum\limits_{n=0}^{\infty}}
\frac{1}{m!n!}\frac{\Gamma(a+m+n)}{\Gamma(a)}\frac{\Gamma(b+m)}{\Gamma
(b)}\frac{\Gamma(b^{\prime}+n)}{\Gamma(b^{\prime})}\frac{\Gamma(c)}%
{\Gamma(c+m+n)}x^{m}y^{n}.
\]
Its numerical implementation in Mathematica is crucial to everything that follows.

\section{Isosceles Triangles \ \ }

Let
\[
T_{\theta}=\left\{  x+iy\in\mathbb{C}:0<x<\cos(\theta/2)\text{, \ }%
y<\tan(\theta/2)x\text{, \ }y>-\tan(\theta/2)x\right\}  ,
\]
the interior of an isosceles triangle with apex angle $0<\theta<\pi$ located
at the origin. \ Two sides of unit length meet there; the third (vertical)
side has $x$-intercept $\cos(\theta/2)$ and length $2$ $\sin(\theta/2)$. \ 

To construct a conformal map $g:\Delta\rightarrow T_{\theta}^{c}$ requires two
steps. \ First, define%
\[
f(z)=%
{\displaystyle\int\limits_{z_{0}}^{z}}
\frac{(\zeta-a_{1})^{\mu_{1}}(\zeta-a_{2})^{\mu_{2}}(\zeta-a_{3})^{\mu_{3}}%
}{\zeta^{2}}d\zeta
\]
on $\Delta$, where $a_{1}$, $a_{2}$, $a_{3}$ are \textit{prevertices} of the
unit circle mapping onto vertices of $T_{\theta}$, and $\pi(1+\mu_{1})$,
$\pi(1+\mu_{2})$, $\pi(1+\mu_{3})$ are exterior angles at the corresponding
vertices. Also, $z_{0}$ is some point of $\Delta$ other than $0$, and the
integral is taken along any curve in $\Delta$ joining $z_{0}$ to $z$ not
passing through $0$ (it does not matter which). Clearly%
\[%
\begin{array}
[c]{ccc}%
\mu_{1}=\dfrac{\pi-\theta}{\pi}, &  & \mu_{2}=\mu_{3}=\dfrac{\pi+\theta}{2\pi
}.
\end{array}
\]
The choice of point $a_{1}$ is arbitrary; here let $a_{1}=-1$. The remaining
two points $a_{2}$, $a_{3}$ must satisfy the constraint \cite{Nh, Ra2}
\[
\mu_{1}/a_{1}+\mu_{2}/a_{2}+\mu_{3}/a_{3}=0
\]
in order that $f(0)=\infty$. Thus
\[%
\begin{array}
[c]{ccc}%
a_{2}=\dfrac{\pi-\theta+2i\sqrt{\pi\theta}}{\pi+\theta}, &  & a_{3}=\dfrac
{\pi-\theta-2i\sqrt{\pi\theta}}{\pi+\theta}%
\end{array}
\]
work for our purposes. We further choose $z_{0}=-1$, so that%
\[
f(z)=%
{\displaystyle\int\limits_{-1}^{z}}
\frac{(\zeta+1)^{1-\theta/\pi}(\zeta^{2}-2\,\tfrac{\pi-\theta}{\pi+\theta
}\,\zeta+1)^{(\pi+\theta)/(2\pi)}}{\zeta^{2}}d\zeta.
\]
The image of $\{a_{1},a_{2},a_{3}\}$ under $f$ evidently lies in the left half
plane -- needing rotation by $\pi$ -- plus rescaling so that the vertical
triangle side has the proper length. \ This second step is achieved by
defining%
\[
g(z)=-\frac{2\sin(\theta/2)}{\operatorname{Im}f(a_{2})-\operatorname{Im}%
f(a_{3})}f(z).
\]
For the scenario $\theta=\pi/2$, it is true that the coefficient
\[
-\frac{\sqrt{2}}{\operatorname{Im}f(a_{2})-\operatorname{Im}f(a_{3}%
)}=0.4756344438799819320567570...=\frac{3^{3/4}}{2^{7/2}\pi^{3/2}}%
\Gamma\left(  \frac{1}{4}\right)  ^{2}=\kappa
\]
to high numerical precision. \ This expression (the outer radius of an
isosceles right triangle) is well-known and is a special case of a more
general formula due to Haegi \cite{PS0, Hg}. \ More on this will be given soon.

Our key result is that the function $f(z)$ possesses an exact representation.
\ Let%
\[%
\begin{array}
[c]{ccc}%
\xi(z)=\dfrac{\left(  \sqrt{\pi}-i\sqrt{\theta}\right)  (z+1)}{2\sqrt{\pi}}, &
& \eta(z)=\dfrac{\left(  \sqrt{\pi}+i\sqrt{\theta}\right)  (z+1)}{2\sqrt{\pi}%
},
\end{array}
\]%
\[
\varphi(z)=\left[  \dfrac{\pi\left(  z-1\right)  ^{2}+\theta(z+1)^{2}}%
{\pi+\theta}\right]  ^{(\pi+\theta)/(2\pi)},
\]%
\[
\psi(z)=\left[  \dfrac{\pi\left(  z-1\right)  ^{2}+\theta(z+1)^{2}}{\pi
+\theta}\right]  ^{(\pi-\theta)/(2\pi)},
\]%
\[
\delta(z)=\left[  -\sqrt{\pi}\left(  z-1\right)  -i\sqrt{\theta}(z+1)\right]
^{(\pi-\theta)/(2\pi)}\left[  -\sqrt{\pi}\left(  z-1\right)  +i\sqrt{\theta
}(z+1)\right]  ^{(\pi-\theta)/(2\pi)}.
\]
Then we have
\begin{align*}
f(z)  &  =(z+1)^{1-\theta/\pi}\left\{  -\frac{\varphi(z)}{z}+\frac
{2^{\theta/\pi}\pi^{(\pi+\theta)/(2\pi)}}{2\pi^{2}+\pi\theta-\theta^{2}}%
\frac{\delta(z)}{\psi(z)}\cdot\right. \\
&  \left[  -2(2\pi-\theta)F_{1}\left(  1-\frac{\theta}{\pi},\frac{\pi-\theta
}{2\pi},\frac{\pi-\theta}{2\pi},2-\frac{\theta}{\pi};\xi(z),\eta(z)\right)
+\right. \\
&  \left.  \left.  (\pi+\theta)(z+1)F_{1}\left(  2-\frac{\theta}{\pi}%
,\frac{\pi-\theta}{2\pi},\frac{\pi-\theta}{2\pi},3-\frac{\theta}{\pi}%
;\eta(z),\xi(z)\right)  \right]  \right\}
\end{align*}
as can be easily proved after-the-fact by differentiation. (Our
before-the-fact technique consisted of examining rational multiples of
$\theta$ in Mathematica, seeking recognizable patterns.) \ This integral
evaluation appears to be new.

Returning to the $\theta=\pi/2$ scenario,
\[
g(z)=\frac{\kappa}{z}+\lambda+O\left(  z\right)
\]
as $z\rightarrow0$, where%
\begin{align*}
\lambda &  =\kappa\frac{2^{5/4}}{3^{3/4}}\left[  2F_{1}\left(  \frac{1}%
{2},\frac{1}{4},\frac{1}{4},\frac{3}{2};\frac{2-i\sqrt{2}}{4},\frac
{2+i\sqrt{2}}{4}\right)  -F_{1}\left(  \frac{3}{2},\frac{1}{4},\frac{1}%
{4},\frac{5}{2};\frac{2+i\sqrt{2}}{4},\frac{2-i\sqrt{2}}{4}\right)  \right] \\
&  =0.5045039334500261012764068...
\end{align*}
is the outer conformal center of $T_{\pi/2}$. \ We wonder whether this
expression for $\lambda$ can be simplified, for example, as a ratio of gamma
or Gauss hypergeometric function values.

As a corollary, let $\tilde{T}=\left\{  x+iy\in\mathbb{C}:x>0\text{,
}y>0\text{, }x+y<1\right\}  $, the initial triangle examined in \cite{Fi1,
Fi2}. \ The outer conformal center of $\tilde{T}$ is simply%
\[
\left(  \frac{1+i}{\sqrt{2}}\right)  \lambda
=(0.3567381524778001406751307...)(1+i)
\]
which is not the same as the inner conformal center $(0.301...)(1+i)$. \ 

We mention finally that the outer conformal center of $T_{\pi/3}$ (an
equilateral triangle) is $1/\sqrt{3}$, that is, it coincides with the centroid
of $T_{\pi/3}$. \ No other scenarios with such recognizable $\lambda$ have
been found.

\section{Haegi's Formula}

An arbitrary triangle with sides $a$, $b$, $c$ and opposite angles
\[
\alpha=\arccos\left(  \frac{b^{2}+c^{2}-a^{2}}{2bc}\right)  ,
\]%
\[
\beta=\arccos\left(  \frac{a^{2}+c^{2}-b^{2}}{2ac}\right)  ,
\]%
\[
\gamma=\arccos\left(  \frac{a^{2}+b^{2}-c^{2}}{2ab}\right)
\]
has area, circumradius and logarithmic capacity given by
\[
A=\sqrt{\frac{a+b+c}{2}\frac{-a+b+c}{2}\frac{a-b+c}{2}\frac{a+b-c}{2}},
\]%
\[
R=\frac{abc}{\sqrt{(a+b+c)(b+c-a)(c+a-b)(a+b-c)}},
\]%
\[
\kappa=\frac{A}{4\pi^{2}q(\alpha/\pi)q(\beta/\pi)q(\gamma/\pi)R}%
\]
where%
\[
q(x)=\frac{1}{\Gamma(x)}\sqrt{\frac{x^{x}}{(1-x)^{1-x}}}.
\]
Under the special circumstances that $a=b=1$ and $c=2$ $\sin(\theta/2)$, we
have%
\[
\kappa(\theta)=\frac{\sqrt{\pi+\theta}}{8\pi^{5/2}}\left(  \frac{\pi+\theta
}{4\theta}\right)  ^{\theta/(2\pi)}\frac{\sin(\theta)^{2}}{\sin(\theta
/2)}\Gamma\left(  \frac{\theta}{\pi}\right)  \Gamma\left(  \frac{\pi-\theta
}{2\pi}\right)  ^{2}%
\]
for the isosceles triangles $T_{\theta}$. \ Over such triangles, the one with
maximal $\kappa$ has $\theta=2.5360873621...$, which seems not to have been
noticed before. Over the family of \textit{all} triangles with fixed $A$, the
one with minimal $\kappa$ is equilateral, as proved by P\'{o}lya \&\ Szeg\"{o}
\cite{PS0, PS3, SZ}. \ If we fix perimeter rather than area, then (to the
contrary) the equilateral triangle provides \textit{maximal} $\kappa$.

\section{Addendum:\ $30^{\circ}$-$60^{\circ}$-$90^{\circ}$ Triangle}

Define $T=\left\{  x+iy\in\mathbb{C}:x>0\text{, }y>0\text{, }\sqrt{3}%
x+y<\sqrt{3}\right\}  $. \ Proceeding as before, we obtain%
\[%
\begin{array}
[c]{lllll}%
\mu_{1}=\dfrac{1}{2}, &  & \mu_{2}=\dfrac{2}{3}, &  & \mu_{3}=\dfrac{5}{6}.
\end{array}
\]
The choice of point $a_{1}$ is arbitrary; here let $a_{1}=1$. From $\mu
_{1}/a_{1}+\mu_{2}/a_{2}+\mu_{3}/a_{3}=0$, we deduce that
\[%
\begin{array}
[c]{ccc}%
a_{2}=i, &  & a_{3}=-\dfrac{3}{5}-\dfrac{4}{5}i
\end{array}
\]
work for our purposes. Choosing $z_{0}=1$, it follows that%
\[
f(z)=%
{\displaystyle\int\limits_{1}^{z}}
\frac{(\zeta-1)^{1/2}(\zeta-i)^{2/3}\left(  \zeta+\left(  \frac{3}{5}+\frac
{4}{5}i\right)  \right)  ^{5/6}}{\zeta^{2}}d\zeta.
\]
An exact representation for $f(z)$ in terms of the Appell $F_{1}$ function is
possible. \ From this, we can verify the outer radius expression \cite{PS0,
Hg} \
\[
\kappa=\frac{5^{5/12}}{2^{10/3}\pi^{2}}\Gamma\left(  \frac{1}{3}\right)
^{3}=0.3779137429709558321024882...
\]
to high numerical precision, but have not yet determined the outer conformal
center of $T$.

\section{Addendum:\ $6$-$9$-$13$ Triangle}

The Schwarz-Christoffel toolbox for Matlab \cite{DT, Dr} makes numerical
computations of a conformal map feasible. \ For the triangle with vertices%
\[
0,\;\;\;6,\;\;\;-\frac{13}{3}+\frac{4\sqrt{35}}{3}i
\]
the following code:%

\[%
\begin{array}
[c]{l}%
\text{\texttt{p = polygon([0 6 -13/3+(4*sqrt(35)/3)*i])}}\\
\text{\texttt{f = extermap(p,scmapopt('Tolerance',1e-18))}}\\
\text{\texttt{p = parameters(f)}}\\
\text{\texttt{format long}}\\
\text{\texttt{p.prevertex}}%
\end{array}
\]
gives%
\[%
\begin{array}
[c]{lllll}%
\mu_{1}=0.659, &  & \mu_{2}=0.207, &  & \mu_{3}=0.132,
\end{array}
\]%
\[%
\begin{array}
[c]{lllll}%
a_{1}=1, &  & a_{2}=0.0163-0.9998i, &  & a_{3}=-0.4069+0.9134i.
\end{array}
\]
Closed-form expressions for these exponents and prevertices are possible yet
cumbersome. \ The same is true for the outer radius $\kappa=3.805336$.
\ Determining the outer conformal center (even approximately) remains open.
\ Figures 1, 2, 3 provide conformal map images of ten evenly-spaced concentric
circles in the disk; orthogonal trajectories are also indicated. \ We leave
the task of exploring whether outer conformal centers belong in Kimberling's
database \cite{Kb} to someone else.%

\begin{figure}[ptb]%
\centering
\includegraphics[
height=3.7204in,
width=3.7204in
]%
{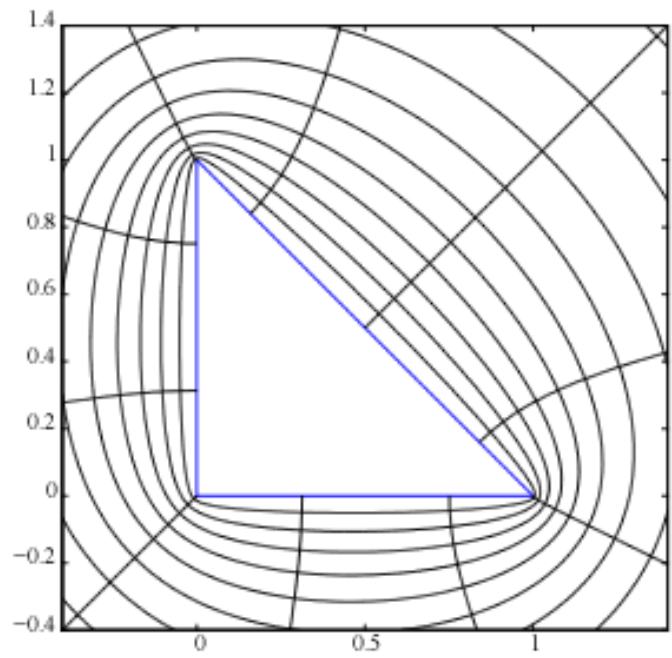}%
\caption{Images of ten concentric circles, center at $0.356+(0.356)i$.}%
\end{figure}
\begin{figure}[ptb]%
\centering
\includegraphics[
height=3.6919in,
width=3.7931in
]%
{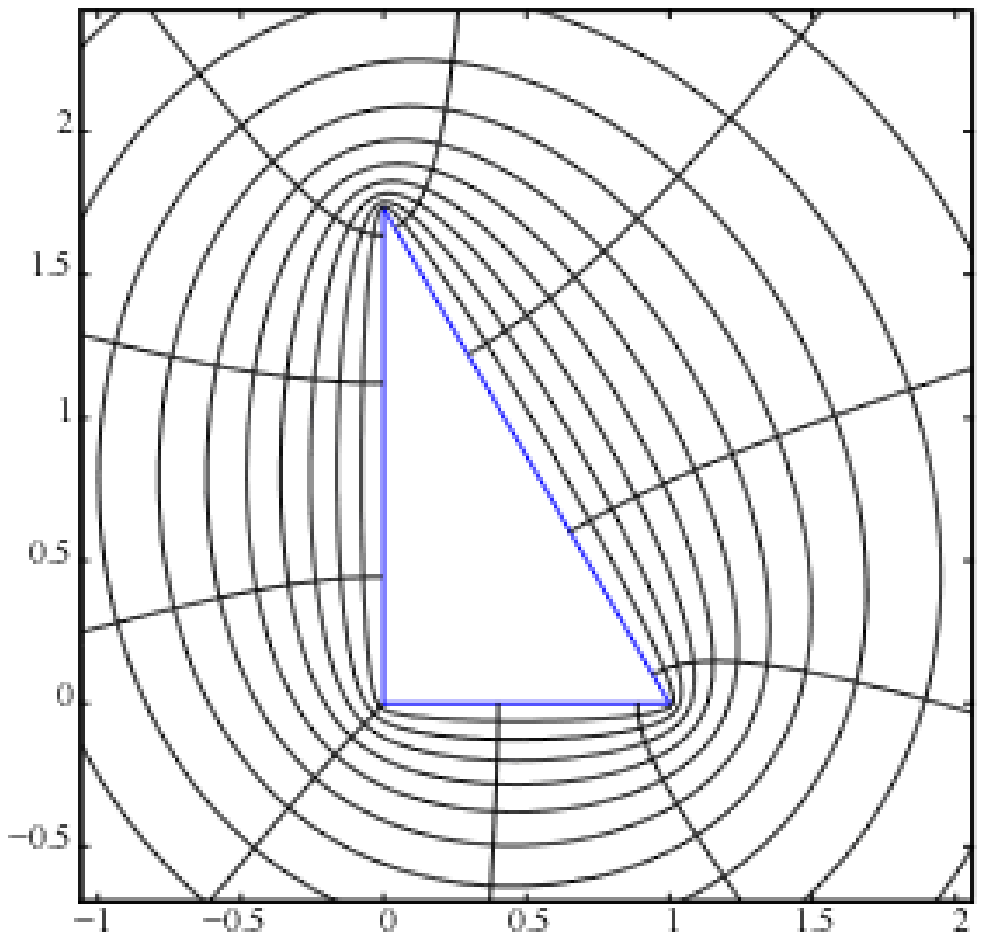}%
\caption{Images of ten concentric circles, center unknown.}%
\end{figure}
\begin{figure}[ptb]%
\centering
\includegraphics[
height=3.781in,
width=3.8458in
]%
{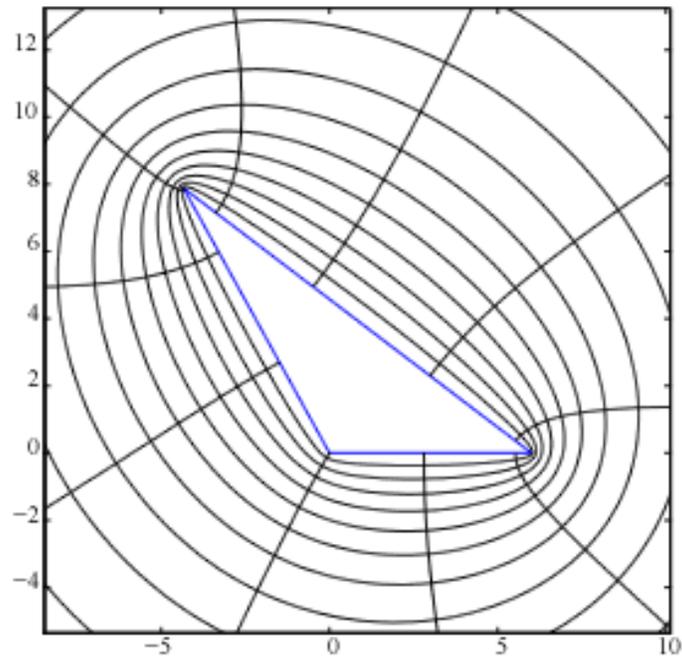}%
\caption{Images of ten concentric circles, center unknown.}%
\end{figure}

\ \ 

\section{Acknowledgements}

I am grateful to Thomas Ransford \cite{Ra1, Ra2} for explaining
Schwarz-Christoffel prevertex calculations for the scenario $\theta=\pi/2$.


\begin{thebibliography}{99}                                                                                               %


\bibitem {Fi1}S. Finch, Least capacity point of triangles, http://arxiv.org/abs/1407.4105.

\bibitem {Kv}J. Kevorkian, \textit{Partial Differential Equations. Analytical
Solution Techniques}, 2$^{\text{nd}}$ ed., Springer-Verlag, 2000, pp.
129--130; MR1728947 (2000i:35001).

\bibitem {PS0}G. P\'{o}lya and\ G. Szeg\"{o}, \textit{Isoperimetric
Inequalities in Mathematical Physics}, Princeton Univ. Press, 1951, pp. 1--3,
158, 254--258, 273--274; MR0043486 (13,270d).

\bibitem {Ha}M. Hantke, \textit{Summen reziproker Eigenwerte}, Ph.D. thesis,
Martin-Luther-Universit\"{a}t Halle-Wittenberg, 2006, http://sundoc.bibliothek.uni-halle.de/diss-online/06/06H308/index.htm.

\bibitem {Pm1}Ch. Pommerenke, On metric properties of complex polynomials,
\textit{Michigan Math. J.} 8 (1961) 97--115; MR0151580 (27 \#1564).

\bibitem {PS1}G. P\'{o}lya and G. Szeg\H{o}, \textit{Problems and Theorems in
Analysis. I, Series, Integral Calculus, Theory of Functions}, Springer-Verlag,
1998, pp. 129--130, 325--326; MR1492447.

\bibitem {PS2}G. P\'{o}lya and G. Szeg\H{o}, \textit{Problems and Theorems in
Analysis. II, Theory of Functions, Zeros, Polynomials, Determinants, Number
Theory, Geometry}, Springer-Verlag, 1998, pp. 23--24, 192--194; MR1492448.

\bibitem {Pm2}Ch. Pommerenke, \textit{Univalent Functions}, Vandenhoeck \&
Ruprecht, 1975, pp. 12--13; MR0507768 (58 \#22526).

\bibitem {BWY}R. Bouffanais, G. D. Weymouth and D. K. P. Yue, Hydrodynamic
object recognition using pressure sensing, \textit{Proc. Royal Soc. London
Ser. A} 467 (2011) 19--38; MR2764670.

\bibitem {Ct}G. Cleanthous, Monotonicity theorems for analytic functions
centered at infinity, \textit{Proc. Amer. Math. Soc.} 142 (2014) 3545--3551.

\bibitem {Lp}S. Liesipohja, \textit{Numerical Methods for Computing
Logarithmic Capacity}, M.Sc. thesis, University of Helsinki, 2014, https://helda.helsinki.fi/handle/10138/44698.

\bibitem {Ra1}T. Ransford, \textit{Potential Theory in the Complex Plane},
Cambridge Univ. Press, 1995, pp. 132--137, 152--160; MR1334766 (96e:31001).

\bibitem {By}W. N. Bailey, \textit{Generalized Hypergeometric Series},
Cambridge Univ. Press, 1935, pp. 73--83.

\bibitem {MSR}J. McDougall, L. Schaubroeck and J. Rolf, \textit{Exploring
Complex Analysis}.\ Ch. 5: \textit{Mappings to Polygonal Domains}, http://www.jimrolf.com/explorationsInComplexVariables.html.

\bibitem {Nh}Z. Nehari,\textit{\ Conformal Mapping}, Dover, 1975, pp.
189--195; MR0377031 (51 \#13206).

\bibitem {Ra2}T. Ransford, Logarithmic capacity for $\tilde{T}$, unpublished
note (2014).

\bibitem {Hg}H. R. Haegi, Extremalprobleme und Ungleichungen konformer
Gebietsgr\"{o}ssen, \textit{Compositio Math.} 8 (1950) 81--111; MR0039811 (12,602b).

\bibitem {Fi2}S. Finch, In limbo: Three triangle centers, http://arxiv.org/abs/1406.0836.

\bibitem {PS3}G. P\'{o}lya and G. Szeg\H{o}, Inequalities for the capacity of
a condenser, \textit{Amer. J. Math.} 67 (1945) 1--32; MR0011871 (6,227e).

\bibitem {SZ}A. Yu. Solynin and V. A. Zalgaller, An isoperimetric inequality
for logarithmic capacity of polygons,\textit{ Annals of Math. }159 (2004)
277--303; MR2052355 (2005a:31002).

\bibitem {DT}T. A. Driscoll and L. N. Trefethen, \textit{Schwarz-Christoffel
Mapping}, Cambridge Univ. Press, 2002, pp. 9--18; MR1908657 (2003e:30012).

\bibitem {Dr}T. A. Driscoll, Schwarz-Christoffel Toolbox for Matlab,
http://www.math.udel.edu/\symbol{126}driscoll/SC/.

\bibitem {Kb}C. Kimberling, Encyclopedia of Triangle Centers, http://faculty.evansville.edu/ck6/encyclopedia/.%

\begin{tabular}
[c]{lll}
& Steven Finch & \\
& Dept. of Statistics & \\
& Harvard University & \\
& Cambridge, MA, USA & \\
& \textit{steven\_finch@harvard.edu} &
\end{tabular}

\end{thebibliography}
\end{document}